\expandafter\ifx\csname amssym.def\endcsname\relax \else\endinput\fi
%
%  Store the catcode of the @ in the csname so that it can be restored later.
\expandafter\edef\csname amssym.def\endcsname{%
       \catcode`\noexpand\@=\the\catcode`\@\space}
%  Set the catcode to 11 for use in private control sequence names.
\catcode`\@=11
%
%  Include all definitions related to the fonts msam, msbm and eufm, so that
%  when this file is used by itself, the results with respect to those fonts
%  are equivalent to what they would have been using AMS-TeX.
%  Most symbols in fonts msam and msbm are defined using \newsymbol;
%  however, a few symbols that replace composites defined in plain must be
%  defined with \mathchardef.

\def\undefine#1{\let#1\undefined}
\def\newsymbol#1#2#3#4#5{\let\next@\relax
 \ifnum#2=\@ne\let\next@\msafam@\else
 \ifnum#2=\tw@\let\next@\msbfam@\fi\fi
 \mathchardef#1="#3\next@#4#5}
\def\mathhexbox@#1#2#3{\relax
 \ifmmode\mathpalette{}{\m@th\mathchar"#1#2#3}%
 \else\leavevmode\hbox{$\m@th\mathchar"#1#2#3$}\fi}
\def\hexnumber@#1{\ifcase#1 0\or 1\or 2\or 3\or 4\or 5\or 6\or 7\or 8\or
 9\or A\or B\or C\or D\or E\or F\fi}

\font\tenmsa=msam10
\font\sevenmsa=msam7
\font\fivemsa=msam5
\newfam\msafam
\textfont\msafam=\tenmsa
\scriptfont\msafam=\sevenmsa
\scriptscriptfont\msafam=\fivemsa
\edef\msafam@{\hexnumber@\msafam}
\mathchardef\dabar@"0\msafam@39
\def\dashrightarrow{\mathrel{\dabar@\dabar@\mathchar"0\msafam@4B}}
\def\dashleftarrow{\mathrel{\mathchar"0\msafam@4C\dabar@\dabar@}}

\def\ulcorner{\delimiter"4\msafam@70\msafam@70 }
\def\urcorner{\delimiter"5\msafam@71\msafam@71 }
\def\llcorner{\delimiter"4\msafam@78\msafam@78 }
\def\lrcorner{\delimiter"5\msafam@79\msafam@79 }
\def\yen{{\mathhexbox@\msafam@55 }}
\def\checkmark{{\mathhexbox@\msafam@58 }}
\def\circledR{{\mathhexbox@\msafam@72 }}
\def\maltese{{\mathhexbox@\msafam@7A }}

\font\tenmsb=msbm10
\font\sevenmsb=msbm7
\font\fivemsb=msbm5
\newfam\msbfam
\textfont\msbfam=\tenmsb
\scriptfont\msbfam=\sevenmsb
\scriptscriptfont\msbfam=\fivemsb
\edef\msbfam@{\hexnumber@\msbfam}

\def\widehat#1{\setbox\z@\hbox{$\m@th#1$}%
 \ifdim\wd\z@>\tw@ em\mathaccent"0\msbfam@5B{#1}%
 \else\mathaccent"0362{#1}\fi}
\def\widetilde#1{\setbox\z@\hbox{$\m@th#1$}%
 \ifdim\wd\z@>\tw@ em\mathaccent"0\msbfam@5D{#1}%
 \else\mathaccent"0365{#1}\fi}
\font\teneufm=eufm10
\font\seveneufm=eufm7
\font\fiveeufm=eufm5
\newfam\eufmfam
\textfont\eufmfam=\teneufm
\scriptfont\eufmfam=\seveneufm
\scriptscriptfont\eufmfam=\fiveeufm

%  Restore the catcode value for @ that was previously saved.
\csname amssym.def\endcsname

\expandafter\ifx\csname pre amssym.tex at\endcsname\relax \else \endinput\fi
%%     Otherwise we store the catcode of the @ in the csname.
\expandafter\chardef\csname pre amssym.tex at\endcsname=\the\catcode`\@
%%     Set the catcode to 11 for use in private control sequence
%%     names.
\catcode`\@=11
%  Most symbols in fonts msam and msbm are defined using \newsymbol.  A few
%  that are delimiters or otherwise require special treatment have already
%  been defined as soon as the fonts were loaded.  Finally, a few symbols
%  that replace composites defined in plain must be undefined first.
\newsymbol\boxdot 1200
\newsymbol\boxplus 1201
\newsymbol\boxtimes 1202
\newsymbol\square 1003
\newsymbol\blacksquare 1004
\newsymbol\centerdot 1205
\newsymbol\lozenge 1006
\newsymbol\blacklozenge 1007
\newsymbol\circlearrowright 1308
\newsymbol\circlearrowleft 1309
\undefine\rightleftharpoons
\newsymbol\rightleftharpoons 130A
\newsymbol\leftrightharpoons 130B
\newsymbol\boxminus 120C
\newsymbol\Vdash 130D
\newsymbol\Vvdash 130E
\newsymbol\vDash 130F
\newsymbol\twoheadrightarrow 1310
\newsymbol\twoheadleftarrow 1311
\newsymbol\leftleftarrows 1312
\newsymbol\rightrightarrows 1313
\newsymbol\upuparrows 1314
\newsymbol\downdownarrows 1315
\newsymbol\upharpoonright 1316
 \let\restriction\upharpoonright
\newsymbol\downharpoonright 1317
\newsymbol\upharpoonleft 1318
\newsymbol\downharpoonleft 1319
\newsymbol\rightarrowtail 131A
\newsymbol\leftarrowtail 131B
\newsymbol\leftrightarrows 131C
\newsymbol\rightleftarrows 131D
\newsymbol\Lsh 131E
\newsymbol\Rsh 131F
\newsymbol\rightsquigarrow 1320
\newsymbol\leftrightsquigarrow 1321
\newsymbol\looparrowleft 1322
\newsymbol\looparrowright 1323
\newsymbol\circeq 1324
\newsymbol\succsim 1325
\newsymbol\gtrsim 1326
\newsymbol\gtrapprox 1327
\newsymbol\multimap 1328
\newsymbol\therefore 1329
\newsymbol\because 132A
\newsymbol\doteqdot 132B
 
\newsymbol\triangleq 132C
\newsymbol\precsim 132D
\newsymbol\lesssim 132E
\newsymbol\lessapprox 132F
\newsymbol\eqslantless 1330
\newsymbol\eqslantgtr 1331
\newsymbol\curlyeqprec 1332
\newsymbol\curlyeqsucc 1333
\newsymbol\preccurlyeq 1334
\newsymbol\leqq 1335
\newsymbol\leqslant 1336
\newsymbol\lessgtr 1337
\newsymbol\backprime 1038
\newsymbol\risingdotseq 133A
\newsymbol\fallingdotseq 133B
\newsymbol\succcurlyeq 133C
\newsymbol\geqq 133D
\newsymbol\geqslant 133E
\newsymbol\gtrless 133F
\newsymbol\sqsubset 1340
\newsymbol\sqsupset 1341
\newsymbol\vartriangleright 1342
\newsymbol\vartriangleleft 1343
\newsymbol\trianglerighteq 1344
\newsymbol\trianglelefteq 1345
\newsymbol\bigstar 1046
\newsymbol\between 1347
\newsymbol\blacktriangledown 1048
\newsymbol\blacktriangleright 1349
\newsymbol\blacktriangleleft 134A
\newsymbol\vartriangle 134D
\newsymbol\blacktriangle 104E
\newsymbol\triangledown 104F
\newsymbol\eqcirc 1350
\newsymbol\lesseqgtr 1351
\newsymbol\gtreqless 1352
\newsymbol\lesseqqgtr 1353
\newsymbol\gtreqqless 1354
\newsymbol\Rrightarrow 1356
\newsymbol\Lleftarrow 1357
\newsymbol\veebar 1259
\newsymbol\barwedge 125A
\newsymbol\doublebarwedge 125B
\undefine\angle
\newsymbol\angle 105C
\newsymbol\measuredangle 105D
\newsymbol\sphericalangle 105E
\newsymbol\varpropto 135F
\newsymbol\smallsmile 1360
\newsymbol\smallfrown 1361
\newsymbol\Subset 1362
\newsymbol\Supset 1363
\newsymbol\Cup 1264
 
\newsymbol\Cap 1265
 
\newsymbol\curlywedge 1266
\newsymbol\curlyvee 1267
\newsymbol\leftthreetimes 1268
\newsymbol\rightthreetimes 1269
\newsymbol\subseteqq 136A
\newsymbol\supseteqq 136B
\newsymbol\bumpeq 136C
\newsymbol\Bumpeq 136D
\newsymbol\lll 136E
 
\newsymbol\ggg 136F
 
\newsymbol\circledS 1073
\newsymbol\pitchfork 1374
\newsymbol\dotplus 1275
\newsymbol\backsim 1376
\newsymbol\backsimeq 1377
\newsymbol\complement 107B
\newsymbol\intercal 127C
\newsymbol\circledcirc 127D
\newsymbol\circledast 127E
\newsymbol\circleddash 127F
\newsymbol\lvertneqq 2300
\newsymbol\gvertneqq 2301
\newsymbol\nleq 2302
\newsymbol\ngeq 2303
\newsymbol\nless 2304
\newsymbol\ngtr 2305
\newsymbol\nprec 2306
\newsymbol\nsucc 2307
\newsymbol\lneqq 2308
\newsymbol\gneqq 2309
\newsymbol\nleqslant 230A
\newsymbol\ngeqslant 230B
\newsymbol\lneq 230C
\newsymbol\gneq 230D
\newsymbol\npreceq 230E
\newsymbol\nsucceq 230F
\newsymbol\precnsim 2310
\newsymbol\succnsim 2311
\newsymbol\lnsim 2312
\newsymbol\gnsim 2313
\newsymbol\nleqq 2314
\newsymbol\ngeqq 2315
\newsymbol\precneqq 2316
\newsymbol\succneqq 2317
\newsymbol\precnapprox 2318
\newsymbol\succnapprox 2319
\newsymbol\lnapprox 231A
\newsymbol\gnapprox 231B
\newsymbol\nsim 231C
\newsymbol\ncong 231D
\newsymbol\diagup 231E
\newsymbol\diagdown 231F
\newsymbol\varsubsetneq 2320
\newsymbol\varsupsetneq 2321
\newsymbol\nsubseteqq 2322
\newsymbol\nsupseteqq 2323
\newsymbol\subsetneqq 2324
\newsymbol\supsetneqq 2325
\newsymbol\varsubsetneqq 2326
\newsymbol\varsupsetneqq 2327
\newsymbol\subsetneq 2328
\newsymbol\supsetneq 2329
\newsymbol\nsubseteq 232A
\newsymbol\nsupseteq 232B
\newsymbol\nparallel 232C
\newsymbol\nmid 232D
\newsymbol\nshortmid 232E
\newsymbol\nshortparallel 232F
\newsymbol\nvdash 2330
\newsymbol\nVdash 2331
\newsymbol\nvDash 2332
\newsymbol\nVDash 2333
\newsymbol\ntrianglerighteq 2334
\newsymbol\ntrianglelefteq 2335
\newsymbol\ntriangleleft 2336
\newsymbol\ntriangleright 2337
\newsymbol\nleftarrow 2338
\newsymbol\nrightarrow 2339
\newsymbol\nLeftarrow 233A
\newsymbol\nRightarrow 233B
\newsymbol\nLeftrightarrow 233C
\newsymbol\nleftrightarrow 233D
\newsymbol\divideontimes 223E
\newsymbol\varnothing 203F
\newsymbol\nexists 2040
\newsymbol\Finv 2060
\newsymbol\Game 2061
\newsymbol\mho 2066
\newsymbol\eth 2067
\newsymbol\eqsim 2368
\newsymbol\beth 2069
\newsymbol\gimel 206A
\newsymbol\daleth 206B
\newsymbol\lessdot 236C
\newsymbol\gtrdot 236D
\newsymbol\ltimes 226E
\newsymbol\rtimes 226F
\newsymbol\shortmid 2370
\newsymbol\shortparallel 2371
\newsymbol\smallsetminus 2272
\newsymbol\thicksim 2373
\newsymbol\thickapprox 2374
\newsymbol\approxeq 2375
\newsymbol\succapprox 2376
\newsymbol\precapprox 2377
\newsymbol\curvearrowleft 2378
\newsymbol\curvearrowright 2379
\newsymbol\digamma 207A
\newsymbol\varkappa 207B
\newsymbol\Bbbk 207C
\newsymbol\hslash 207D
\undefine\hbar
\newsymbol\hbar 207E
\newsymbol\backepsilon 237F
%  Restore the catcode value for @ that was previously saved.
\catcode`\@=\csname pre amssym.tex at\endcsname

%\endinput
%\input mssymb      % these must be input at the BEGINNING of this file; some
% MAXIMUM LINE LENGTH IS 80 CHARACTERS
%
%\font\eurofont=times at 10pt% contains Û symbol (as option-shift-2)
%
\font\fivebi=cmmib5
\font\fivebsy=cmbsy5
\font\sixrm=cmr6
\font\sixi=cmmi6
\font\sixbf=cmbx6
\font\sixsy=cmsy6
\font\sixmsa=msam5 at 6pt% This is needed since Bluesky's virtual fonts for
\font\sixmsb=msbm5 at 6pt% sizes 6.8.9 all come out too small (same below)
\font\sevenbi=cmmib7
\font\sevenbsy=cmbsy7
\font\eightrm=cmr8
\font\eightsl=cmsl8
\font\eightit=cmti8
\font\eighti=cmmi8
\font\eightbf=cmbx8
\font\eightsy=cmsy8
\font\eightmsa=msam7 at 8pt
\font\eightmsb=msbm7 at 8pt
\font\ninerm=cmr9
\font\ninesl=cmsl9
\font\nineit=cmti9
\font\ninei=cmmi9
\font\ninebf=cmbx9
\font\ninebi=cmmib10 scaled 900
\font\ninesy=cmsy9
\font\ninebsy=cmbsy10 scaled 900
\font\ninemsa=msam10 at 9pt
\font\ninemsb=msbm10 at 9pt
\font\tenbit=cmbxti10
\font\tenbsl=cmbxsl10
\font\tenbi=cmmib10
\font\tenbsy=cmbsy10
\font\twelvebf=cmbx12
\font\twelvebi=cmmib10 scaled 1200
\font\twelvebsy=cmbsy10 at 12pt

\let\sc=\sevenrm            % SMALL CAPS (in tenpoint)
\def\eightpoint{%
     %\font\eurofont=times at 8pt% contains Û symbol (as option-shift-2)
     \def\rm{\fam0\eightrm}%         see p.415
     \textfont0=\eightrm \scriptfont0=\sixrm \scriptscriptfont0=\fiverm
     \textfont1=\eighti \scriptfont1=\sixi \scriptscriptfont1=\fivei
     \textfont2=\eightsy \scriptfont2=\sixsy \scriptscriptfont2=\fivesy
     \textfont3=\tenex \scriptfont3=\tenex \scriptscriptfont3=\tenex
     \textfont\itfam=\eightit \def\it{\fam\itfam\eightit}%
     \textfont\slfam=\eightsl \def\sl{\fam\slfam\eightsl}%
     \textfont\bffam=\eightbf \scriptfont\bffam=\sixbf
     \scriptscriptfont\bffam=\fivebf \def\bf{\fam\bffam\eightbf}%
     \textfont\msbfam=\eightmsb \textfont\msafam=\eightmsa
     \scriptfont\msafam=\sixmsa \scriptfont\msbfam=\sixmsb
     \scriptscriptfont\msafam=\fivemsa \scriptscriptfont\msbfam=\fivemsb
      \skewchar\eighti='177 \skewchar\sixi='177
      \skewchar\eightsy='60 \skewchar\sixsy='60
     \normalbaselineskip=10pt%(normally 9pt; must be =height+depth in next line)
     \setbox\strutbox=\hbox{\vrule height7pt depth3pt width0pt}%(normally 7&2pt)
     \let\sc=\sixrm \let\big=\eightbig \normalbaselines\rm}
\def\ninepoint{%
     %\font\eurofont=times at 9pt% contains Û symbol (as option-shift-2)
     \def\rm{\fam0\ninerm}%         see p.415
     \textfont0=\ninerm \scriptfont0=\sixrm \scriptscriptfont0=\fiverm
     \textfont1=\ninei \scriptfont1=\sixi \scriptscriptfont1=\fivei
     \textfont2=\ninesy \scriptfont2=\sixsy \scriptscriptfont2=\fivesy
     \textfont3=\tenex \scriptfont3=\tenex \scriptscriptfont3=\tenex
     \textfont\itfam=\nineit \def\it{\fam\itfam\nineit}%
     \textfont\slfam=\ninesl \def\sl{\fam\slfam\ninesl}%
     \textfont\bffam=\ninebf \scriptfont\bffam=\sixbf
     \scriptscriptfont\bffam=\fivebf \def\bf{\fam\bffam\ninebf}%
     \textfont\msbfam=\ninemsb \textfont\msafam=\ninemsa
     \scriptfont\msafam=\sixmsa \scriptfont\msbfam=\sixmsb
     \scriptscriptfont\msafam=\fivemsa \scriptscriptfont\msbfam=\fivemsb
      \skewchar\ninei='177 \skewchar\sixi='177
      \skewchar\ninesy='60 \skewchar\sixsy='60
     \normalbaselineskip=11pt
     \setbox\strutbox=\hbox{\vrule height8pt depth3pt width0pt}%
     \let\sc=\sevenrm \let\big=\ninebig \normalbaselines\rm}
\def\tenpoint{%
     %\font\eurofont=times at 10pt% contains Û symbol (as option-shift-2)
     \def\rm{\fam0\tenrm}%         see p.415
     \textfont0=\tenrm \scriptfont0=\sevenrm \scriptscriptfont0=\fiverm
     \textfont1=\teni \scriptfont1=\seveni \scriptscriptfont1=\fivei
     \textfont2=\tensy \scriptfont2=\sevensy \scriptscriptfont2=\fivesy
     \textfont3=\tenex \scriptfont3=\tenex \scriptscriptfont3=\tenex
     \textfont\itfam=\tenit \def\it{\fam\itfam\tenit}%
     \textfont\slfam=\tensl \def\sl{\fam\slfam\tensl}%
     \textfont\bffam=\tenbf \scriptfont\bffam=\sevenbf
     \scriptscriptfont\bffam=\fivebf \def\bf{\fam\bffam\tenbf}%
     \textfont\msbfam=\tenmsb \textfont\msafam=\tenmsa
     \scriptfont\msafam=\sevenmsa \scriptfont\msbfam=\sevenmsb
     \scriptscriptfont\msafam=\fivemsa \scriptscriptfont\msbfam=\fivemsb
     \normalbaselineskip=12pt
%     \normalbaselineskip=13pt% (normally 12pt; here 13 = 9 + 4 of next line)
%     \setbox\strutbox=\hbox{\vrule height9pt depth4pt width0pt}%
     \let\sc=\sevenrm \let\big=\tenbig \normalbaselines\rm}
\catcode`@=11        % This allows the use of `@' in the next line (p.344)
\def\tenbig#1{{\hbox{$\left#1\vbox to8.5pt{}\right.\n@space$}}}
\def\ninebig#1{{\hbox{$\textfont0=\tenrm\textfont2=\tensy
      \left#1\vbox to7.25pt{}\right.\n@space$}}}
\def\eightbig#1{{\hbox{$\textfont0=\ninerm\textfont2=\ninesy
      \left#1\vbox to6.5pt{}\right.\n@space$}}}
\catcode`@=12        % This restores the `inhibiting' catcode of `@'.
\def\bold{%
     \textfont0=\tenbf \scriptfont0=\sevenbf \scriptscriptfont0=\fivebf
     \textfont1=\tenbi \scriptfont1=\sevenbi \scriptscriptfont1=\fivebi
     \textfont2=\tenbsy \scriptfont2=\sevenbsy \scriptscriptfont2=\fivebsy
       \textfont\itfam=\tenbit \def\it{\fam\itfam\tenbit}%
       \textfont\slfam=\tenbsl \def\sl{\fam\slfam\tenbsl}%
       \textfont\bffam=\tenbf \scriptfont\bffam=\sevenbf
       \textfont\msbfam=\tenmsb \textfont\msafam=\tenmsa
   \fam0\tenbf}
\def\bigbold{%
     \textfont0=\twelvebf \scriptfont0=\ninebf \scriptscriptfont0=\sevenbf
     \textfont1=\twelvebi \scriptfont1=\ninebi \scriptscriptfont1=\sevenbi
     \textfont2=\twelvebsy \scriptfont2=\ninebsy \scriptscriptfont2=\sevenbsy
     \fam0\twelvebf}
%
%
% MAXIMUM LINE LENGTH BELOW IS 80 CHARACTERS
%
% Instructions on how to use \item, and how to surround items by space,
% are given in paperformat.
%
\let\plainitem=\item% \item itself is redefined in paperformat.
\let\plainitemitem=\itemitem% \itemitem itself is redefined in paperformat.
%

%
% Û is option-shift-2 on the Mac
%
                %  needs mssymb
              %  in old files, \nat may be used instead

%
%
%                            SKIPS AND BREAKS
%
\def\g{\hskip.17em\relax}               %  breakable thin space
\def\th{\thinspace}                     %  non-breakable thin space
\def\nl{\hfil\break}
\newskip\Bigskipamount
   \Bigskipamount=2\baselineskip plus.5\baselineskip minus.3\baselineskip
\def\Bigbreak{\removelastskip\vskip0pt plus .1\vsize\penalty-1000
              \vskip0pt plus-.1\vsize\vskip\Bigskipamount}
% The following can be used after a display to remove the belowdisplayskip, e.g.
% when the display finishes a proof and is just followed by \endproof. If text
% follows, e.g. just a word like "and" leading to another display, one may wish
% to say \noskip\smallskip\noindent.
%use after display followed by short line
\def\Nobreak$$#1$${\postdisplaypenalty=10000$$#1$$\postdisplaypenalty=0}
%
%
%                             ABBREVIATIONS
%
     % \H will redefined below.

\let\doublebar=\| 
\def\|{\!\!\restriction\!\!}

  \let\sub=\sube

\def\supe{\supseteq}

\def\sm{\smallsetminus}
\def\es{\emptyset}

\def\T{T} \def\H{I}% italics, like K^n, P^n etc. Cf. I FOR "Inflated"
\def\wrt{with respect to}

\def\:{\colon}
\def\minor{\preccurlyeq} 
\def\Minor{\succcurlyeq}

\def\slt{\mathrel{\hbox{$\minor$\kern-.6em\lower.33ex\hbox{${}_s\;$}}}}
\def\sgt{\mathrel{\mathchoice                        %("simplicial minors")
   {\hbox{$\Minor$\kern-.5em\lower.3ex\hbox{${}_s$}}}
   {\hbox{$\Minor$\kern-.5em\lower.3ex\hbox{${}_s$}}}
   {\hbox{$\scriptstyle\Minor\kern-.43em\lower.28ex\hbox{$\scriptstyle{}_s$}$}}
 {\hbox{$\scriptstyle\Minor\kern-.43em\lower.28ex\hbox{$\scriptstyle{}_s$}$}} }}    

\def\ucl(#1){\lfloor #1 \rfloor}% up-closure
\def\dcl(#1){\lceil #1 \rceil}% down-closure
%
%

%

% Example: '\interior P' puts a circle over the P.
% (the "\relax" is important: otherwise a capital letter A-F instead of P
% is interpreted as part of the number 7017,
% giving a "wrong math code complaint!!)
%
% The following commands can be used to 'specify' a relation (= #1) by putting
% something (= #2) underneath or above.
% Examples:\specrel<T, \specrel\sim G, \specrel={(1)}, \Specrel\Rightarrow?.
%
\def\specrel#1#2{\mathrel{\mathop{\kern0pt #1}\limits_{#2}}}
\def\Specrel#1#2{\mathrel{\mathop{\kern0pt #1}\limits^{#2}}}
%
% Next a version of \specrel for use in aligned equations;
% here the specification text is put in an \hbox of width 0,
% so as not to interfere with the alignment.
%
\def\alignspecrel#1#2{\mathrel{\mathop{\kern0pt #1}\limits_{\hbox
   to0pt{\hss$\scriptstyle#2$\hss}}}}
\def\alignSpecrel#1#2{\mathrel{\mathop{\kern0pt #1}\limits^{\hbox
   to0pt{\hss$\scriptstyle#2$\hss}}}}
\def\invlim{\specrel\lim{\raise 2pt\hbox{$\longleftarrow$}}}% inverse limit, projective limit
\def\proof{\removelastskip\penalty55\medskip\noindent{\bf Proof. }}
\def\noproof{{\unskip\nobreak\hfill\penalty50\hskip2em\hbox{}\nobreak\hfill%
       $\square$\parfillskip=0pt\finalhyphendemerits=0\par}\goodbreak}
\def\endproof{\noproof\bigskip}
% Syntax example for the following: \looseproof{Zweiter Beweis von Satz \xxxA}
\def\looseproof#1{\bigbreak\noindent {{\bf #1.}}}%space follows in input

\newcount\refno
\def\ref#1#2\par{{\plainitem{[??]}#2\smallskip}}
\newtoks\thingtowrite %USED AGAIN LATER FOR INDEX
\long\def\writerefnumber#1{%
    \thingtowrite={#1}%
    \immediate\write\refnumbersfile{\the\thingtowrite}%
    }
\newwrite\refnumbersfile
\def\makerefnumbers{\immediate\openout\refnumbersfile=RefNumbers%
  \refno=0 \writerefnumber{\refno=0}
  \def\ref##1##2\par{\global\advance\refno by 1
    \writerefnumber{\global\advance\refno by 1 \newcounter##1 ##1=\the\refno}%
       % \newcounter replaces \newcount, which is forbidden inside a def.
       % When using the auto-generated file, say "\let\newcounter=\newcount"
       % before reading in that file.
    \plainitem{[\the\refno]}##2\smallskip% automatic numbering, ignoring ##1
    }%
  }
\def\autorefnumbers{\refno=0
  \def\ref##1##2\par{\advance\refno by 1\plainitem{[\the\refno]}##2\smallskip}
  }
\def\userefnumbers{\refno=0
  \def\ref##1##2\par{\advance\refno by 1\plainitem{[\the##1]}##2\smallskip}
  }
% TO USE, FIRST RUN WITH \makerefnumbers active (and no file named "RefNumbers" present
% in the directory to which such a file is to be written; it cannot overwrite).
% THEN DISABLE "\makerefnumbers", and make the following line active:
% \userefnumbers\let\newcounter=\newcount\input RefNumbers
%
%
\def\proclaimwithname #1. (#2) #3\par{{\bigbreak
  \clubpenalty=10000\noindent{\bf#1.\enspace}(#2)\nl
  {\sl #3}\par\bigbreak}}
\def\proposition (#1) #2\par{{\setbox0\hbox{(#1)\enspace}\bigbreak
   \sl\hangindent\the\wd0 \noindent\hskip\the\wd0
   \llap{\box0}\ignorespaces#2\par\bigbreak}}
\def\subsection #1\par{\vskip 3\medskipamount minus \smallskipamount\leftline{\bold #1}
        \penalty10000\smallskip\noindent}
      \def\section #1\par{\Bigbreak\centerline{\bf #1} % TO BE PHASED OUT
              \penalty10000\bigskip\noindent}
%
%         The following versions of the beginsection macro are for
%         section headings immediately followed by \proclaim:
\def\beginpsection #1\par{\Bigbreak\centerline{\bold #1}
        \penalty10000\bigskip}
\def\psubsection #1\par{\bigbreak\leftline{\bold #1}\penalty10000\bigskip}
%
%         The following macro positions its argument flush right,
%         like the \endproof box or an equation number.

%
%  The following three items usually need an \enditem. They are for use inside
%  \proclaim or another \item. Note that these items don't require their own
%  \par at the end, but they have no finishing \smallskip. (There seems to be
%  no way around this: we can't delimit by \par, since in the input this would
%  end the \proclaim.) So, if the item is followed by some text within the same
% \proclaim (say), the \enditem should (and can) be replaced with \smallskip.
%        
\def\pitem#1{\smallskip\advance\parindent by 3mm
             \plainitem{\rm(#1)}\advance\parindent by-3mm}
\def\pitemitem#1{\smallskip\advance\parindent by 3mm
             \plainitemitem{\rm(#1)}\advance\parindent by-3mm}
\def\varitemitem#1{{\setbox0\hbox{\hskip\parindent#1\enskip}
           \smallbreak\hangindent\the\wd0 \noindent\hskip\the\wd0
           \llap{#1\enskip}\ignorespaces}}
%
%        Note that the following variations of \item must end with \par.
%
\newdimen\newparindent
\def\iitem#1#2\par{\newparindent=\parindent \advance\newparindent by 3mm
           \smallbreak \hangindent\newparindent \noindent\hskip\newparindent
           \llap{{\rm #1}\enspace}\ignorespaces#2\par\smallbreak}
\def\iitemitem#1#2\par{\newparindent=\parindent \advance\newparindent by 3mm
           \smallbreak \hangindent2\newparindent \noindent\hskip2\newparindent
           \llap{{\rm #1}\enspace}\ignorespaces#2\par\smallbreak}
\def\varitem#1#2\par{{\setbox0\hbox{{\rm #1}\enspace}
           \smallbreak \hangindent\the\wd0 \noindent\hskip\the\wd0
           \llap{{\rm #1}\enspace}\ignorespaces#2\par\smallbreak}}
\def\enditem{\par}
% For use in Exercises and Hints
%
\def\Textindent#1{\par \advance\parindent by 3mm
                  \textindent{{\rm #1}} \advance\parindent by -3mm}
\def\indentedline#1{\advance\hsize by -\parindent \line{#1}
                   \advance\hsize by \parindent}
\def\iindentedline#1{\advance\parindent by 3mm
                     \advance\hsize by -\parindent
                     \line{#1}
                     \advance\hsize by \parindent
                     \advance\parindent by -3mm}
\newdimen\margin   % needed for macros \textdisplay & \ltextdisplay
%  The following macro takes 3 arguments, #1 and #3 in math-mode,
%  #2 as plain text. It displays #1\quad centered w.r.t. the whole \line
%  (if #2 leaves enough space), adds \quad#2 to the right of #1\quad, and puts
%  #3 flush right. The arguments must be separated by &'s. The argu-
%  ments themselves may be empty, but there must be 2 ampersands.
%  Examples: $$\textdisplay x=y &for all $x\in X$& (1')$$
%            $$\textdisplay x=y &for all $x\in X$&$$
\def\textdisplay#1&#2&#3$${\margin=\hsize
          \setbox1=\hbox{$\displaystyle#1\quad$}%
          \setbox2=\hbox{\quad#2\qquad$#3$}%
                     \advance\margin by-\wd1
                     \divide\margin by 2
   \ifdim\wd2 < \margin
      \box1\rlap{\quad#2}\eqno#3$$%
   \else
      \line{\qquad\hfil \box1\quad #2 \qquad $#3$}$$%
   \fi}
%
%  The following macro is the `\leqno' version of \textdisplay; argument
%  #1 is the \leqno, #2 is the formula to be displayed, and #3 is the
%  text following the formula.
\def\ltextdisplay#1&#2&#3$${\margin=\hsize
           \setbox2=\hbox{$\displaystyle#2\quad$}
           \setbox3=\hbox{\quad#3\qquad}
                     \advance\margin by-\wd2
                     \divide\margin by 2
   \ifdim\wd3 < \margin
      \line{$#1$\hfil\box2\hbox to \margin{\box3\hfil}}$$%
   \else
      \line{$#1$\qquad\hfil\box2\quad #3\qquad} $$%
   \fi}
%
%   The next macro displays and centres #1, typically a paragraph of text.
%   #2, typically an `eqno', is set flush right, vertically centred, and
%   has to be in math mode.
%   Example: \textno This is a ... lot of text&(3')\par
\def\textno#1&#2\par{%
   \margin=\hsize
   \advance\margin by -4\parindent
          \setbox1=\hbox{\sl#1}%
   \ifdim\wd1 < \margin
      $$\box1\eqno#2$$\endgraf%
   \else
      \bigbreak
      \line{\indent$\vcenter{\advance\hsize by -3\parindent
      \sl\noindent#1}\hfil#2$}%
      \bigbreak
   \fi}
%
%   The same with a left eqno:
\def\textlno#1&#2\par{%
   \margin=\hsize
   \advance\margin by -4\parindent
          \setbox1=\hbox{\sl#1}%
   \ifdim\wd1 < \margin
      $$\box1\leqno#2$$%
   \else
      \bigbreak
      \line{$#2\hfil\vcenter{\advance\hsize by -3\parindent
          \sl\noindent#1}\hskip\parindent$}%
      \bigbreak
   \fi}
%
%
%                        MARGINAL HACKS ETC.
%
\newcount\commentno
\def\COMMENT#1{$^{<\the\commentno>}$%
     \vadjust{\vbox to 0pt{\vss\vskip-8pt\rightline{%
     \rlap{\hbox{\hskip7mm \vbox{\pretolerance=-1
     \doublehyphendemerits=0 \finalhyphendemerits=0
     \hsize40mm\tolerance=10000\eightpoint
     \lineskip=0pt\lineskiplimit=0pt
     \rightskip=0pt plus16mm\baselineskip8pt\noindent
     \hskip0pt       %(without this, the first word is never hyphenated!)
     {$\langle$\the\commentno. #1$\rangle$}\endgraf}}}}\vss}}%
     \global\advance\commentno by1}%
\def\writecommentsasfootnotes{%
 \def\COMMENT{\global\advance\commentno by1\footnote{$^{<\the\commentno>}$}}%
 }
\def\nocomments{\def\COMMENT##1{}}
%
%
%   The \? macro puts
%   the argument #1 in the left margin. Examples: \??, \?{What nonsense!}.
\def\?#1{\vadjust{\vbox to 0pt{\vss\vskip-8pt\leftline{%
     \llap{\hbox{\vbox{\pretolerance=-1
     \doublehyphendemerits=0\finalhyphendemerits=0
     \hsize16truemm\tolerance=10000\eightpoint
     \lineskip=0pt\lineskiplimit=0pt
     \rightskip=0pt plus16truemm\baselineskip8pt\noindent
     \hskip0pt        %(without this, the first word is never hyphenated!)
     #1\endgraf}\hskip7truemm}}}\vss}}}
\def\d{}% DISABLES WHAT'S JUST BEEN DEFINED ABOVE
%
% The following \ds macro is the "silent" version of \d: it writes its argument
% in the margin (or index file), just as \d does, but not into the current text.
% Note the different syntax: no delimiting blank in input.
%
%  CHANGE: I got fed up with marginal reminders to define things, so this
%  feature is disabled right away by the line "\def\ds#1{}" below. Note that
%  \makeindex still works, because it defines \d and \ds anew.
%
%\def\ds#1{\ifmmode
%     \vadjust{\vbox to 0pt{\vss\vskip-8pt\leftline{%
%     \llap{\hbox{\vbox{\pretolerance=-1
%    \doublehyphendemerits=0\finalhyphendemerits=0
%     \hsize16truemm\tolerance=10000\eightpoint
%     \lineskip=0pt\lineskiplimit=0pt
%     \rightskip=0pt plus16truemm\baselineskip8pt\noindent
%     define $#1$!\endgraf}\hskip7truemm}}}\vss}}%
%   \else
%     \vadjust{\vbox to 0pt{\vss\vskip-8pt\leftline{%
%     \llap{\hbox{\vbox{\pretolerance=-1
%     \doublehyphendemerits=0\finalhyphendemerits=0
%     \hsize16truemm\tolerance=10000\eightpoint
%     \lineskip=0pt\lineskiplimit=0pt
%     \rightskip=0pt plus16truemm\baselineskip8pt\noindent
%     \hskip0pt        %(without this, the first word is never hyphenated!)
%     define ``#1''!\endgraf}\hskip7truemm}}}\vss}}%
%   \fi}
%
\def\ds#1{}% DISABLES WHAT'S JUST BEEN DEFINED ABOVE
%
% The following is a device (due to CET1) which allows to \write something to
% a file without expanding all the tokens completely. (This would result in
% the use of lots of "@"s, which make the file untexable. Another way around the
% problem would be to precede every argument of \write by "\catcode`@=11" (and
% to set it back to 12 after the argument). That would make the file texable,
% but since it would still be difficult to read, the method below is better.
% NOTE: This works well for a Remarks file (say), but NOT for an index file. 
% The reason is that in order to get the page numbers of entries from the first
% paragraph of a page right, one has to say "\write" rather than
% "\immediate\write". But in the def of \indexwrite one has to say "\immediate",
% since otherwise the token gets overwritten and rather
% than n index words from a page the index will
% contain the last index word of that page n times. But if the index words are
% written to file immediately while their page numbers get queued, then the two
% get separated from one another.
%
%\newtoks\thingtowrite % NOW EARLIER
\long\def\indexwrite#1{%
    \thingtowrite={#1}%
    \immediate\write\index{\the\thingtowrite}%
    }
%
%  The following macro \makeindex redefines \d and \ds.
%  It suppresses the appearance of \d's arguments in the margin
%  and writes them to a file called "index" instead.
%
\newwrite\index
\def\makeindex{\immediate\openout\index=index%
   \immediate\write\index{\catcode`@=11}%
   \def\d##1 {\ifmmode
     \write\index{$##1$, }%
     \write\index{\the\count0}\write\index{}% blank line for para
   \else
     \write\index{{##1}, }%
     \write\index{\the\count0}\write\index{}% blank line for para
   \fi {##1} }% It's important to put the text AFTER the def: otherwise
              % the blank is not merged with possible blanks following in
              % the input file, which may result in an additional blank line
      \def\ds##1{\ifmmode
     \write\index{$##1$, }%
     \write\index{\the\count0}\write\index{}% blank line for para
   \else
     \write\index{##1, }%
     \write\index{\the\count0}\write\index{}% blank line for para
   \fi}}
%
%   The \m macro puts its (direct!) argument in the right margin,
%   leaving it also in the text (in italics). 
%   Syntax example: this \m{word\/} is defined here
%   and will be shown in the margin; similarly the set $X$ in $A =: \m X$.
%   Intended use: in preprints, at places where something is first defined
%   (to enhance readability of the paper, in the absence of an index).
%   To disable, say "\disablems" at beginning of file.
\newdimen\gap% gap between text and hack
\gap=3truemm
\newdimen\hackwidth
\hackwidth=15truemm
\def\disablems{\def\mos##1{\strut}}% NOT simply "{}": when \m is active
 % it puts a strut in the line, so this should be here also when \m is disabled
 % to avoid a change in vertical space (and hence possibly in page breaks)
\def\m#1{\ds{#1}\mo{#1}}% \ds writes to index file, \mo in margin and text
\def\mo#1{\ifmmode {#1}\else {\it#1}\fi\mos{#1}}
\def\ms#1{\ds{#1}\mos{#1}}% \ds writes to index file, \mos in margin
\def\mos#1{\ifmmode
     \strut\vadjust{\vbox to 0pt{\vss\kern-11pt\leftline{%
     \llap{\hbox{\vbox{\pretolerance=-1
     \doublehyphendemerits=0\finalhyphendemerits=0
     \hsize\hackwidth\tolerance=10000\eightpoint
     \lineskip=0pt\lineskiplimit=0pt
     \rightskip=0pt plus\hsize\baselineskip8pt\noindent
     $#1$\strut\endgraf}\hskip\gap }}}\vss}}%
   \else
     \strut\vadjust{\vbox to 0pt{\vss\kern-11pt\leftline{%
     \llap{\hbox{\vbox{\pretolerance=-1
     \doublehyphendemerits=0\finalhyphendemerits=0
     \hsize\hackwidth\tolerance=10000\eightpoint
     \lineskip=0pt\lineskiplimit=0pt
     \rightskip=0pt plus\hsize\baselineskip8pt\noindent
     \hskip0pt    %(without this, the first word is never hyphenated!)
     {\sl#1}\strut\endgraf}\hskip\gap }}}\vss}}%
   \fi}%
\newcount\remarkno
\def\REMARK#1{{\footnote{${}^{\the\remarkno}$}{{#1}}%
   \global\advance\remarkno by1}}
\def\noremarks{\def\REMARK##1{}}
%
%
%
%                              PICTURES
%
\def\picture #1 by #2 (#3){
  \vbox to #2{
          \vfill
          \special{picture #3}
          \hrule width #1 height 0pt depth 0pt
           }}
\newdimen\topfiguremargin
   \topfiguremargin=0pt                                  % default
\newdimen\bottomfiguremargin
   \bottomfiguremargin=\medskipamount                    % default
\newdimen\normalpictureheight
\normalpictureheight=40mm
%   The following macro uses TeXtures pictures; these MUST be named
%   Fig.1, Fig.2.5 etc. (without a space after the '.'), as in the macro call
%itself. The width (#2) and height (#3) should have the original values
%of the TeXtures picture, to facilitate proper scaling. The heightfactor (#4)
%is divided by 1000 and used to scale the \normalpictureheight. Thus, height-
%factor 500 (2000) results in a picture of half (twice) the normalpicureheight.
%   Syntax example for a figure at the standard \normalpictureheight:
%\Fig.2 (538pt by 536pt; heightfactor: 1000; caption: This is a short caption)
\def\Fig.#1 (#2by#3; heightfactor:#4; caption:#5) {{%
   \dimen2=\normalpictureheight
   \dimen0=#2                          % computing width
      \divide\dimen2 by 1000
      \multiply\dimen2 by#4              % \dimen2 := intended pictureheight
   \count2=\dimen2                  % computing scalefactor
      \dimen1=#3                             % \dimen1 := actual pictureheight
   \count1=\dimen1
   \divide\count1 by 1000
   \divide\count2 by \count1          % \count2 := scalefactor (times 1000)
%      \message{scalefactor in Fig.#1 is \the\count2}%
   \divide\dimen0 by 1000
   \multiply\dimen0 by \count2      % \dimen0 := width
         \dimen1=\hsize
         \advance\dimen1 by -\dimen0
         \divide\dimen1 by 2               % \dimen1 := margin
   \midinsert
   \vbox to \topfiguremargin{\vfil}
   \noindent\hskip\dimen1
   \picture\dimen0 by \dimen2  (Fig.#1 scaled \the\count2)%
   \vskip\bottomfiguremargin                     % beginning caption
      \ninepoint
      \parindent=.1\hsize\narrower\narrower
      \setbox0\hbox{#5}
      \ifdim\wd0 < .6\hsize
           \centerline{F{\sc IGURE} #1.\hskip1em#5}
       \else
           \plainitem{F{\sc IGURE} #1. }#5\par
       \fi
   \vskip0pt\endinsert}}
%
%The following \textpicture macro is for inserting pictures in the current
%text line. The width (#2) and height (#3) should have the original values of
%the TeXtures picture, to facilitate proper scaling. #4 offers the opportunity
%to reserve vertical space for the picture by saying "height15pt depth10pt" or
%so; this will be the height of the \vbox containing the picture (default=0).
%#6 is the amount by which the bottom left corner is placed below the baseline
%(poss.neg.), #5 is the horizontal extension of the picture. Syntax examples:
%   \textpicture flower(538pt by 536pt; width5em lower-5pt)
%   \textpicture bug(538pt by 536pt; height20pt depth20pt width5em lower20pt)
\def\textpicture #1(#2by#3; #4width#5lower#6){{%
  % computing scalefactor
      \dimen0=#5\count2=\dimen0                    % desired width
      \dimen0=#2\count1=\dimen0                    % actual width
   \divide\count1 by 1000
   \divide\count2 by \count1                 % \count2 is now = scalefactor
  %\count3=11820\divide\count3 by \count2
  %\message{The vertices of #1 should have width \the\count3}
   \hbox{\vrule #4width0pt\vbox to 0pt{\vss\vskip#6%
      \special{picture #1 scaled \the\count2}\hrule width#5 height0pt\vss}}}}
%
%
%The following \figure macro builds on the \epsf macro and includes an EPS file.
% THE FOLLOWING LINE HAS TO OCCUR AT THE START OF THE MAIN FILE TO BE TEXED:
% \input epsf.def
%
% #1 is just a figure number to be used in the caption.
% #2 is the caption itself.
% #3 is the file name of the figure; this must an EPS file
%    (or PS with bounding box).
% #4 is for scaling, but disabled now (see below)
%
% Syntax example: \figure 3. The funny graph $G$ (Funny.graph.eps; 800)
%
\def\figure #1. #2 (#3; #4) {{%
   \def\bigskip{\par\ifdim\lastskip<\bigskipamount\removelastskip
                                              % eg. for abb after \endproof
      \vskip\bigskipamount\fi}% takes effect inside the def(!) of midinsert
   \midinsert\vskip\topfiguremargin
   \dimen0=\normalpictureheight
      \divide\dimen0 by 1000
      \multiply\dimen0 by#4        % \dimen0 := intended pictureheight
   \centerline{\epsfbox{#3.eps}}%                  good placing - use this!
   \vskip\bottomfiguremargin                     % beginning caption
      \ninepoint
      \parindent=.1\hsize\narrower\narrower
      \setbox0\hbox{#2}
      \ifdim\wd0 < .6\hsize
           \centerline{F{\sc IGURE} #1.\hskip1em#2}
       \else
           \plainitem{F{\sc IGURE} #1. }#2\par
       \fi
  \endinsert}}
%
% The following is for notes for talks to be hidden on transparancies.
% Say \hide{bla} to make "bla" disappear (but it will take up space);
% \showhidden will show all hidded text in grey (eg for paper notes).
%

%
%                               REFERENCES
%
\def\Abh#1 {{\sl Abh.\g Math.\g Sem.\g Univ.\g Hamburg\penalty100\ \bf#1\ }}
\def\AMASH#1 {{\sl Acta Math.\g Acad.\g Sci.\g Hung.\penalty100\ \bf#1\ }}
\def\Advances#1 {{\sl Adv.\g Math.\penalty100\ \bf#1\ }}
\def\Annals#1 {{\sl Ann.\g Math.\penalty100\ \bf#1\ }}
\def\AnnComb#1 {{\sl Ann.\g Comb.\penalty100\ \bf#1\ }}
\def\AMM#1 {{\sl Amer.\g Math.\g Monthly\penalty100\ \bf#1\ }}
\def\Archiv#1 {{\sl Arch.}\g {\sl Math.\penalty100\ \bf#1\ }}
\def\ArsComb#1 {{\sl Ars Comb.\penalty100\ \bf#1\ }}% FROMERLY \AC
\def\CJM#1 {{\sl Can.\g J.\th Math.\penalty100\ \bf#1\ }}
\def\Comb#1 {{\sl Com\-bi\-na\-to\-ri\-ca\penalty100\ \bf#1\ }}
\def\CPC#1 {{\sl Comb.\g Probab.\g Comput.\penalty100\ \bf#1\ }}
\def\Crelle#1 {{\sl J.}\th {\sl Reine Angew.}\g
    {\sl Math.\penalty100\ \bf#1\ }}
\def\DM#1 {{\sl Discrete Math.\penalty100\ \bf#1\ }}
\def\DAM#1 {{\sl Discrete Appl.\g Math.\penalty100\ \bf#1\ }}
\def\EJC#1 {{\sl Eur.}\g{\sl J.}\g{\sl Comb.\penalty100\ \bf#1\ }}
\def\EJ#1 {{\sl Electronic.}\g{\sl J.}\g{\sl Comb.\penalty100\ \bf#1\ }}
\def\GC#1 {{\sl Graphs Comb.\penalty100\ \bf#1\ }}
\def\IJ#1 {{\sl Isr.\g J.\th Math.\penalty100\ \bf#1\ }}
\def\Inv#1 {{\sl In\-vent.\g math.\penalty100\ \bf#1\ }}
\def\JAlg#1 {{\sl J.}\th {\sl Algorithms\penalty100\ \bf#1\ }}
\def\JCTA#1 {{\sl J.}\th {\sl Comb.}\g {\sl Theory~A\penalty100\ \bf#1\ }}
\def\JCTB#1 {{\sl J.}\th {\sl Comb.}\g {\sl Theory~B\penalty100\ \bf#1\ }}
\def\JGT#1 {{\sl J.}\th {\sl Graph Theory\penalty100\ \bf#1\ }}
\def\BLMS#1 {{\sl Bull.\g Lond.\g Math.\g Soc.\penalty100\ \bf#1\ }}
\def\JLMS#1 {{\sl J.\g Lond.\g Math.\g Soc.\penalty100\ \bf#1\ }}
\def\PLMS#1 {{\sl Proc.\g Lond.\g Math.\g Soc.\penalty100\ \bf#1\ }}
\def\Order#1 {{\sl Order\ \bf#1\ }}
\def\Random#1 {{\sl Random Struct.\g Alg.\penalty100\ \bf#1\ }}
\def\MA#1 {{\sl Math.}\g {\sl Ann.\penalty100\ \bf#1\ }}
\def\MN#1 {{\sl Math.}\g {\sl Nachr.\penalty100\ \bf#1\ }}
\def\MPCPS#1 {{\sl Math.\g Proc.\g Camb.\g Phil.\g Soc.\penalty100\ \bf#1\ }}
\def\MS#1 {{\sl Math.}\g {\sl Scand.\penalty100\ \bf#1\ }}
\def\MZ#1 {{\sl Math.}\g {\sl Zeit.\penalty100\ \bf#1\ }}
\def\BAMS#1 {{\sl Bull.\th Amer.\g Math.\g Soc.\penalty100\ \bf#1\ }}
\def\JAMS#1 {{\sl J.\th Amer.\g Math.\g Soc.\penalty100\ \bf#1\ }}
\def\MAMS#1 {{\sl Mem.\g Amer.\g Math.\g Soc.\penalty100\ \bf#1\ }}
\def\PAMS#1 {{\sl Proc.\g Amer.\g Math.\g Soc.\penalty100\ \bf#1\ }}
\def\SIAM#1 {{\sl SIAM J.\g Discrete Math.\penalty100\ \bf#1\ }}
\def\SLNM#1 {{\sl Springer Lecture Notes in Mathematics\penalty100\ \bf#1\ }}
\def\TAMS#1 {{\sl Trans.\g Amer.\g Math.\g Soc.\penalty100\ \bf#1\ }}
\def\TCSA#1 {{\sl Theor.\g Comput.\g Sci.~A\penalty100\ \bf#1\ }}
%
%
%
%
%
%                     ONLY FOR TeXtures without mssymb.tex:
%
%\def\N{{\rm \rlap I{\kern.18em}N}}   
%  \catcode`@=11        % This allows the use of `@' in the next line (p.344)
%\def\not#1{\mathrel{\mathpalette\c@ncel#1}} % use (only) with Imagewriter.
%  \catcode`@=12        % This restores the `inhibiting' catcode of `@'.
%\def\subsetneqq{\mathrel %{\mathchoice
%   {\vcenter{
%      \hbox{\lower6pt\hbox{$\scriptstyle\subset$}}
%      \hbox{\raise3pt\hbox{$\flatneq$}}}} }
%   \def\flatneq{\rlap {$\scriptstyle =$} {\kern1.5pt} {\scriptscriptstyle /}}
%\def\square{\Square53}            % See def. of \Square below.
%\def\Square#1#2{{\vbox{\hrule height.#2pt
%       \hbox{\vrule width.#2pt height#1pt \kern#1pt
%          \vrule width.#2pt}
%       \hrule height.#2pt}}}
%\def\nexists{\hbox{\rm\raise1.2pt\rlap/$\exists$}}
%
%
%
%
%
%\def\language=#1{}% For use with Textures versions < 1.3
%
%
%                  MODIFICATIONS TO PLAIN TeX
%
%\catcode`[=\active \catcode`]=\active
%  \def[{\thinspace\lbrack\thinspace}
%  \def]{\thinspace\rbrack}
%\def\{{\lbrace\thinspace}
%\def\}{\thinspace\rbrace}
%
\bigskipamount=1\baselineskip plus.3\baselineskip minus.3\baselineskip
\medskipamount=\bigskipamount\divide\medskipamount by 2
\smallskipamount=\medskipamount\divide\smallskipamount by 2 % (p.349)
\medmuskip = 3mu plus 2mu minus 1mu
\thickmuskip = 6mu plus 4mu minus 2mu % (for spacing in formulae; pp.168/170)
\def\smallbreak{\par \ifdim\lastskip<\smallskipamount
   \removelastskip \penalty-100 \smallskip \fi}
\def\medbreak{\par \ifdim\lastskip<\medskipamount
   \removelastskip \penalty-250 \medskip \fi}
\def\bigbreak{\par \ifdim\lastskip<\bigskipamount
   \removelastskip \penalty-500 \bigskip \fi}
\catcode`@=11        % This allows the use of `@' in the next line (p.344)
  \def\raggedbottom{\topskip10pt plus20pt \r@ggedbottomtrue} % The amount of
%                      permitted raggedness is controlled by the `plus' item.
\catcode`@=12        % This restores the `inhibiting' catcode of `@'.
\def\ge{\geqslant}% \geq remains available for the default version of \ge.
\def\le{\leqslant}% \leq remains available for the default version of \le.
\let\elt=\in
\def\in{\mathrel{\mathchoice
   {\raise .7pt \hbox{$\scriptstyle\elt$}}
   {\raise .7pt \hbox{$\scriptstyle\elt$}}
   {\raise .5pt \hbox{$\hskip .5pt\scriptscriptstyle\elt\hskip .5pt$}}
   {\raise.35pt \hbox{$\scriptscriptstyle\elt$}} }}
\let\hasaselt=\owns
\def\owns{\mathrel{\mathchoice
   {\raise .7pt \hbox{$\scriptstyle\hasaselt$}}
   {\raise .7pt \hbox{$\scriptstyle\hasaselt$}}
   {\raise .5pt \hbox{$\hskip .5pt\scriptscriptstyle\hasaselt\hskip .5pt$}}
   {\raise.35pt \hbox{$\scriptscriptstyle\hasaselt$}} }}
\let\exis=\exists
   \def\exists{\exis\>}
\let\nexis=\nexists
   \def\nexists{\nexis\>}
                            % To be phased out
\let\foral=\forall
   \def\forall{\foral\>}
\let\Rightarro=\Rightarrow
   \def\Rightarrow{\>\Rightarro\>}
\let\mi=\min
   \def\min{\mi\>}
\let\ma=\max
   \def\max{\ma\>}
\let\su=\sup
   \def\sup{\su\>}
\let\inff=\inf
   \def\inf{\inff\>}
\mathchardef\to="2221   % = \rightarrow, but of `binop' type (p.154)
\def\proclaim #1.#2 #3\par{\bigbreak
   \noindent{\bf#1.}#2\enspace{\sl#3}\par\bigbreak}
% The second argument above is optional and intended for references. Note that
% it is terminated by a space, which must therefore be present unsuppressed in
% input. Syntax examples: "\proclaim Thm {1.3}. This is the theorem.\par" or
% "\proclaim Thm \Euler.\five{} This is Euler's theorem.\par", where \Euler
%  expands to {1.3} and \five to \th [5], say. (Note that \proclaim Thm 1.3.
% would treat the 3. as argument #2 (incorrectly) and fail to set it in bold.)
% Or directly: "\proclaim Theorem \xxxVTop.~[\the\ref] blabla" (note the ~).
%
\newskip\sectionheadlineskipamount
\sectionheadlineskipamount=8pt plus 2pt minus 1pt
\def\beginsection #1\par{\Bigbreak\centerline{\bold #1}
        \penalty10000\vskip\sectionheadlineskipamount\noindent}
\let\ffootnote=\footnote
\def\footnote#1#2{\ffootnote{#1}{\eightpoint#2\vskip-12pt}}
%                  (The \vskip inserts an implicit \par, which has two
%                    effects: first, the desired effect of wrapping
%                    up the last paragraph of the footnote giving it the
%                    correct linespacing (that of \eightpoint), secondly
%                    the undesired effect of starting a new paragraph with
%                    a strut. To counteract the arising blank vertical
%                    space, the skip is chosen negative.)
%
% Note that the following redefinition of \item finishes with a smallbreak.
% Since two smallbreaks in sequence result in only a single smallbreak, this
% gives a smallskip at the beginning, between any two items, and at the end.
% if additional space is desired before and after a series of items, say
% \medskip before and \par\smallskip (not \smallbreak) after the series of
% items.
%
\def\item#1#2\par{\parindent=10mm\smallbreak\hang\indent
                  \llap{{\rm #1}\enspace}\ignorespaces#2\par\smallbreak
                  \parindent=7mm}
\def\itemitem#1#2\par{\parindent=10mm\smallbreak
                  \indent\hangindent2\parindent\indent
                  \llap{{\rm #1}\enspace}\ignorespaces#2\par\smallbreak
                  \parindent=7mm}
%
%
%                     FORMAT PARAMETERS
%
\pretolerance=0 %This prevents line breaks in maths formulas - I don't know why.
\tolerance=2000
%\fontdimen2\tenrm=3.8pt
%\fontdimen2\tensl=3.8pt
%\fontdimen2\tenit=4pt
\baselineskip=13pt                 %(DEFAULT IS 12pt)
\vsize=200mm                   % (used to be 240truemm; changed 06/98)
\hsize=120mm                   % (used to be 140truemm; changed 06/98)
\hoffset=9mm                   % (used to be 9truemm; changed 06/98)
\parindent=7mm
\relpenalty=2000 \binoppenalty=5000  % DISCOURAGES BREAKS IN FORMULAS
\hyphenpenalty=100
\abovedisplayskip=12pt plus3pt minus4pt
\belowdisplayskip=12pt plus3pt minus4pt    % (p.348)
%
%   Unfortunately, TeX is unable to adjust the skip following a display to the
%   length of the line **below** the display, in the way in which it chooses
%  between \abovedisplayskip and \abovedisplayshortskip depending on the length
%   of the line above it. Thus, such adjustment has to be done by hand: setting
%
\belowdisplayshortskip=9pt plus3pt minus3pt
    % (formerly 12pt, like \belowdisplayskip)
%
%   reduces the effect of the standard 'short' version, while
%
%                  \def\noskip{\vskip-\lastskip\noindent}
%
%   (which is contained in macros.tex) removes any belowdisplay skip
%   (long or short) altogether.
%   One may want to say \smallskip\noindent just after it.
%
   %\language=\Germanlanguage %Commented out to work with (non-bilingual) Plain, for ArXiv
% DON'T SAY "\German" HERE; THE READING-IN OF GERMAN.TEX WHICH THIS
% CAUSES CREATES ERRORS E.G. WITH LATER USE OF THE " CHARACTER
% IN THE \specialS USED FOR MAKING PDF LINKS, BECAUSE " BECOMES
% ACTIVE AND IS DEFINED IN AN UNDESIRED WAY.
% ***HOWEVER: if, for some reason, it becomes necessary to say \German
% here, it is possible to toggle the "-character between active and normal,
% using the commands \mdqon and \mdqoff defined in German.tex.
%
 % UNFORTUNATELY, NO UMLAUTS CAN BE USED IN \HYPHENATION 
 \hyphenation{Baum-ord-nung Baum-ord-nun-gen End-ecke End-ecken kur-zen
Kur-zen Graphen-ei-gen-schaft Graphen-ei-gen-schaften he-raus he-raus-ar-bei-ten
he-raus-zu-ar-bei-ten Schnitt-raum}%
 %\English %Commented out to work with (non-bilingual) Plain, for ArXiv
 \hyphenation{ac-cess-ible ana-log-ous ana-log-ous-ly ana-lyze ana-lyse
ana-ly-sis answer answers aver-age bundle bundles Buch-ge-sell-schaft col-our
col-ours col-oured col-our-ing col-our-ings con-struct-ible con-struct-ive
con-struct-ive-ly co-rol-lary Co-rol-lary des-cend des-cend-ing Deut-sche
end-li-cher de-fi-ni-tion de-fi-ni-tions De-fi-ni-tion equi-val-ent
equi-val-ence Euler-ian exist-ence every Gra-phen Hamil-ton-ian homeo-mor-phic
homeo-mor-phism homeo-mor-phisms hy-po-thesis hy-po-theses in-ac-cess-ible
ir-regu-lar ir-regu-lar-ity method methods modi-fi-ca-tion mono-chro-matic par-ticu-lar
pro-po-si-tion pro-po-si-tions Pro-po-si-tion regu-lar regu-lar-ity regu-lar-ly
sig-ni-fi-cant sig-ni-fi-cant-ly sig-ni-fi-cance to-po-lo-gical to-po-lo-gical-ly
un-at-tached un-end-li-cher using Using Wis-sen-schaft-li-che}
%
%\autorefnumbers
%\makerefnumbers
\userefnumbers\let\newcounter=\newcount\refno =0
\global \advance \refno by 1 \newcounter \refMetric  \refMetric =\the \refno 
\global \advance \refno by 1 \newcounter \refTopSurvey  \refTopSurvey =\the \refno 
\global \advance \refno by 1 \newcounter \refBook  \refBook =\the \refno 
\global \advance \refno by 1 \newcounter \refBGC  \refBGC =\the \refno 
\global \advance \refno by 1 \newcounter \refNST  \refNST =\the \refno 
\global \advance \refno by 1 \newcounter \refHalin  \refHalin =\the \refno 
\global \advance \refno by 1 \newcounter \refJung  \refJung =\the \refno 

\hbox{}\vskip2pt\centerline{\bigbold A simple existence criterion for normal spanning trees}
\smallskip\centerline{\bigbold in infinite graphs}
\vskip 4mm
\centerline{Reinhard Diestel}
\disablems
\noremarks
\nocomments

\def\NST{normal spanning tree}
\def\TKA{$\T K_{\aleph_0}$}

\def\xxxCharacterizations{1}
\def\xxxMain{2}
\def\xxxAAproperties{3}
\def\xxxAAfat{4}
\def\xxxIAAfat{5}
\def\xxxATfat{6}
\def\xxxIATfat{7}

\def\secResult{1}
\def\secProof{2}

\bigskip\medskip
{\narrower\narrower\ninepoint\noindent
   Halin proved in 1978 that there exists a normal spanning tree in every connected graph $G$ that satisfies the following two conditions: (i)~$G$~contains no subdivision of a `fat'~$K_{\aleph_0}$, one in which every edge has been replaced by uncountably many parallel edges; and (ii)~$G$ has no $K_{\aleph_0}$ subgraph. We show that the second condition is unnecessary.\par}

\beginsection Introduction

A spanning tree of an infinite graph is {\it normal\/} if the endvertices of any chord are comparable in the tree order defined by some arbitrarily chosen root. (In finite graphs, these are their `depth-first search' trees; see~[\the\refBook] for precise definitions.) Normal spanning trees are perhaps the most important single structural tool for analysing an infinite graph~-- see [\the\refBGC] for a typical example, and the exercises in [\the\refBook, Chapter~8] for many more~-- but they do not always exist. The question of which graphs have \NST s thus is an important question.

All countable connected graphs have \NST s~[\the\refBook]. But not all connected graphs do. For example, if $T$ is a \NST\ of~$G$ and $G$ is complete, then $T$ defines a chain on its vertex set. Hence $T$ must be a single path or ray, and $G$ is countable.

For connected graphs of arbitrary order, there are three characterizations of the graphs that admit a \NST:

\proclaim Theorem \xxxCharacterizations.
The following statements are equivalent for connected graphs~$G$.
 \pitem{i} $G$ has a \NST{\/\rm;}
 \pitem{ii} $V(G)$ is a countable union of dispersed sets {\rm ({\sl Jung}~[\the\refJung, \the\refTopSurvey]);}
 \pitem{iii} $|G|$ is metrizable {\rm [\the\refMetric];}
 \pitem{iv} $G$ contains neither an $(\aleph_0,\aleph_1)$-graph nor an Aronszajn-tree graph as a minor {\rm [\the\refNST]}.\enditem

\noindent
Here, a set of vertices in $G$ is {\it dispersed\/} if every ray can be separated from it by some finite set of vertices. (The levels of a \NST\ are dispersed; see~[\the\refBook].) The dispersed vertex sets in a graph~$G$ are precisely those that are closed in the topological space~$|G|$ of~(iii), which consists of $G$ and its ends~[\the\refMetric]. The space $|G|$ will not concern us in this note, so we refer to~[\the\refMetric] for the definition of the topology on~$|G|$. But we shall use the equivalence of (i) and~(iv) in our proof, and the forbidden minors mentioned in~(iv) will be defined in Section~\secProof.

Despite the variety in Theorem~\xxxCharacterizations, it can still be hard in practice to decide whether a given graph has a \NST.%
   \footnote{$^1$}{In particular, the two types of graph mentioned in~(iv) are not completely understood; see~[\the\refNST] for the~-- quite intriguing~-- problem of how to properly understand (or meaningfully classify) the $(\aleph_0,\aleph_1)$-graphs.}
   In most applications, none of these characterizations is used, but a simpler sufficient condition due to Halin. This condition, however, is much stronger, and hence does not always hold even if a \NST\ exists. It is the purpose of this note to show that this condition can be considerably weakened.

\beginsection \secResult. The result

Halin's~[\the\refHalin] most-used sufficient condition for the existence of a \NST\ in a connected graph is that it does not contain a $\T K_{\aleph_0}$. This is usually easier to check than the conditions in Theorem~\xxxCharacterizations, but it is also quite a strong assumption. However, Halin~[\the\refHalin] also proved that this assumption can be replaced by the conjunction of two independent much weaker assumptions:\medskip

\item{$\bullet$} $G$ contains no {\it fat $\T K_{\aleph_0}$}: a subdivision of the multigraph obtained from a $K_{\aleph_0}$ by replacing every edge with $\aleph_1$ parallel edges;

\item{$\bullet$} $G$ contains no $K_{\aleph_0}$ (as a subgraph).

\smallskip\noindent
We shall prove that the second condition is unnecessary:

\proclaim Theorem \xxxMain.
Every connected graph not containing a fat $\T K_{\aleph_0}$ has a \NST.

We remark that all the graphs we consider are simple, including our fat $\T K_{\aleph_0}$s. When we say, without specifying any graph relation, that a graph $G$ {\it contains\/} another graph~$H$, we mean that $H$ is isomorphic to a subgraph of~$G$. Any other undefined terms can be found in~[\the\refBook].

\beginsection \secProof. The proof

Our proof of Theorem~\xxxMain\ will be based on the equivalence (i)$\leftrightarrow$(iv) in Theorem~\xxxCharacterizations, so let us recall from~[\the\refNST] the terms involved here.

An {\it Aronszajn tree\/} is a poset $(T,\le)$ with the following properties:%
   \footnote{$^2$}{Unlike the perhaps better known Suslin trees~-- Aronszajn trees in which even every anti\-chain must be countable~-- Aronszajn trees can be shown to exist without any set-theoretic assumptions in addition to~ZFC.}\medskip
  \item{$\bullet$} $T$ that has a least element, its {\it root\/};
  \item{$\bullet$} the down-closure of every point in $T$ is well-ordered;
  \item{$\bullet$} $T$ is uncountable, but all chains and all levels in $T$ are  countable.

\smallskip\noindent
Here, the {\it down-closure\/} $\dcl(t)$ of a point $t\in T$ is the set $\{\,x\mid x\le t\,\}$; its {\it up-closure\/} is the set $\ucl(t) := \{\,y\mid t\le y\,\}$.\vadjust{\penalty-500}
   More generally, if $x<y$ we say that $x$ lies {\it below\/} $y$ and $y$ {\it above\/}~$x$. The {\it height\/} of a point $t\in T$ is the order type of the chain $\dcl(t)\sm\{t\}$, and the {\it levels\/} of $T$ are its maximal subsets of points of equal height.

An {\it Aronszajn-tree graph\/} or {\it AT-graph\/}, is a graph~$G$ on whose vertex set there exists an Aronszajn tree $T$ such that

\medskip
\item{$\bullet$} the endvertices of every edge of $G$ are comparable in~$T$;
\item{$\bullet$} for all $x<y$, the vertex $y$ has a neighbour $x'$ such that $x\le x'<y$.

\smallskip\noindent
The second condition says that each vertex is joined cofinally to the vertices below it. The idea behind this is that if we were to construct any order tree $T$ on $V(G)$ satifying the first condition, a tree satisfying also the second condition would be one that minimizes the level of each vertex.%
   \COMMENT{}

Note that {\it intervals\/} in~$T$, sets of the form $\{\,t\mid x\le t < y\,\}$ for some given points $x<y$, span connected subgraphs in~$G$. This is because every ${t > x}$ has a neighbour $t'$ with $x\le t' < t$, by the second condition, and hence the interval contains for each of its elements~$t$ the vertices of a $t$--$x$ path in~$G$. Similarly, $G$~itself is connected, because every vertex can be linked to the unique root of~$T$.\looseness=-1

An {\it $(\aleph_0,\aleph_1)$-graph with bipartition $(A,B)$\/} is a bipartite graph with vertex classes $A$ of size~$\aleph_0$ and $B$ of size~$\aleph_1$ such that every vertex in $B$ has infinite degree.

Replacing the vertices $x$ of a graph $X$ with disjoint connected graphs~$H_x$, and the edges $xy$ of $X$ with non-empty sets of $H_x$--$\,H_y$ edges, yields a graph that we shall call an~$\H X$ (for `inflated~$X$'). More formally, a graph {\it $H$ is an~$\m{\H X}$} if its vertex set admits a partition
${\{\,V_x\mid x\in V(X)\,\}}$ into connected subsets~$V_x$ such that distinct vertices $x,y\in X$ are adjacent in $X$ if and only if $H$ contains a $V_x$--$V_y$ edge. The sets $V_x$ are the {\it branch sets\ms{branch sets} of the~$\H X$\/}. Thus, $X$~arises from $H$ by contracting the subgraphs~$H_x$, without deleting any vertices or edges (other than loops or parallel edges arising in the contraction). A~graph $X$ is a {\it minor\/} of a graph $G$ if $G$ contains an~$\H X$ as a subgraph. See~[\the\refBook] for more details.

For our proof of Theorem~\xxxMain\ from Theorem~\xxxCharacterizations\ (i)$\leftrightarrow$(iv) it suffices to show the following:

\textno
Every $\H X$, where $X$ is either an $(\aleph_0,\aleph_1)$-graph or an AT-graph, contains a fat $\T K_{\aleph_0}$ (as a subgraph). &(*)

\noindent
The rest of this section is devoted to the proof of~$(*)$.

\proclaim Lemma \xxxAAproperties.
Let $X$ be an $(\aleph_0,\aleph_1)$-graph, with bipartition $(A,B)$ say.
   \pitem{i} $X$ has an $(\aleph_0,\aleph_1)$-subgraph $X'$ with bipartition into $A'\sub A$ and ${B'\sub B}$ such that every vertex in $A'$ has uncountable degree in~$X'$.
   \pitem{ii} For every finite set $F\sub A$ and every uncountable set $U\sub B$, there exists a vertex $a\in A\sm F$ that has uncountably many neighbours in~$U$.
 \enditem

\proof
(i) Delete from $X$ all the vertices in~$A$ that have only countable degree, together with their neighbours in~$B$. Since this removes only countably many vertices from~$B$, the remaining set~$B'\sub B$ is still uncountable.\vadjust{\penalty-2000} Every $b'\in B'$ has all its $X$-neighbours in the set $A'$ of the vertices in~$A$ that we did not delete, as otherwise $b'$ would have been deleted too. Thus, $b'$~still has infinite degree in the subgraph $X'$ of $X$ induced by $A'$ and~$B'$. In particular, $A'$~is still infinite, and $X'$ is the desired $(\aleph_0,\aleph_1)$-subgraph of~$X$.

(ii) If there is no vertex $a\in A\sm F$ as claimed, then each vertex $a\in A\sm F$ has only countably many neighbours in~$U$. As $A\sm F$ is countable, this means that $U\sm N(A\sm F)\ne\es$. But every vertex in this set has all its neighbours in~$F$, and thus has finite degree. This contradicts our assumption that $X$~is an $(\aleph_0,\aleph_1)$-graph.
   \endproof

\proclaim Lemma \xxxAAfat.
Let $X$ be an $(\aleph_0,\aleph_1)$-graph with bipartition $(A,B)$. Let $A'\sub A$ be infinite and such that for every two vertices $a,a'$ in $A'$ there is some uncountable set $B(a,a')$ of common neighbours of $a$ and~$a'$%
   \COMMENT{}
   in~$B$. Then $A'$ is the set of branch vertices of a fat~\TKA\ in $X$ whose subdivided edges all have the form $aba'$ with $b\in B(a,a')$.

\proof
We have to find a total of $\aleph_0^2\cdot\aleph_1 = \aleph_1$ independent paths in $X$ between vertices in~$A'$. Let us enumerate these desired paths as $(P_\alpha)_{\alpha < \omega_1}$; it is then easy to find them recursively on~$\alpha$, keeping them independent. 
   \COMMENT{}
   \endproof

\proclaim Lemma \xxxIAAfat.
Every $\H X$, where $X$ is an $(\aleph_0,\aleph_1)$-graph, contains a fat $\T K_{\aleph_0}$.

\proof
Let $H$ be an $\H X$ for an $(\aleph_0,\aleph_1)$-graph $X$ with bipartition $(A,B)$, with branch sets $V_x$ for vertices $x\in X$. Replacing $X$ with an appropriate $(\aleph_0,\aleph_1)$-subgraph $Y$ (and $H$ with the corresponding~$\H Y\sub H$) if necessary, we may assume by Lemma~\xxxAAproperties\th(i) that every vertex in~$A$ has uncountable degree in~$X$. We shall find our desired fat \TKA\ in~$H$ as follows.

We construct, inductively, an infinite set $A' = \{a_0,a_1,\dots\}\sub A$ such that, for each $a_i\in A'$, there is an uncountable subdivided star $S(a_i)\sub H[V_{a_i}]$ whose leaves send edges of~$H$ to the branch sets of distinct vertices $b\in B$. The sets $B_i$ of these~$b$%
   \COMMENT{}
   will be nested as $B_0\supe B_1\supe\dots$. We shall then apply Lemma~\xxxAAfat\ to find a fat \TKA\ in~$X$, and translate this to the desired fat \TKA\ in~$H$.

Pick $a_0\in A$ arbitrarily. For each of the uncountably many neighbours $b$ of $a_0$ in $B$ we can find a vertex $v_b\in V_b$ that sends an edge of~$H$ to~$V_{a_0}$. For every~$b$, pick one neighbour $u_b$ of $v_b$ in~$V_{a_0}$. Consider a minimal connected subgraph $H_0$ of~$H[V_{a_0}]$ containing all these vertices~$u_b$, and add to it all the edges~$u_b v_b$ to obtain the\vadjust{\penalty-200} graph~$T = T(a_0)$. By the minimality of~$H_0$,
   \textno $T$~is a tree in which every edge lies on a path between two vertices of the form~$v_b$. & (1)%
   \COMMENT{}

Since there are uncountably many~$b$%
   \COMMENT{}
   and their $v_b$ are distinct, $T$~is uncountable and hence has a vertex $s_0$ of uncountable degree. For every edge $e$ of $T$ at~$s_0$ pick a path in $T$ from $s_0$ through $e$ to some~$v_b$; this is possible by~(1). Let $S(a_0)$ be the union of all these paths. Then $S(a_0)$ is an uncountable subdivided star with centre~$s_0$ all whose non-leaves lie in~$V_{a_0}$ and whose leaves lie in the branch sets $V_b$ of distinct vertices $b\in B$. Let $B_0\sub B$ be the (uncountable) set of these~$b$, and rename the vertices $v_b$ with $b\in B_0$ as~$v_b^0$.

Assume now that, for some $n\ge 1$, we have picked distinct vertices $a_0,\dots,a_{n-1}$ from~$A$ and defined uncountable subsets $B_0\supe \dots\supe B_{n-1}$ of~$B$ so that each $a_i$ is adjacent in $X$ to every vertex in~$B_i$. By Lemma~\xxxAAproperties\th (ii) there exists an $a_n\in A\sm \{a_0,\dots,a_{n-1}\}$ which, in~$X$, has uncountably many neighbours in~$B_{n-1}$. As before, we can find an uncountable subdivided star $S(a_n)$ in~$H$ whose centre $s_n$ and any other non-leaves lie in $V_{a_n}$ and whose leaves $v_b^n$ lie in the branch sets $V_b$ of (uncountably many) distinct vertices $b\in B_{n-1}$. We let $B_n$ be the set of those~$b$. Then $B_n$ is an uncountable subset of~$B_{n-1}$, and $a_n$~is adjacent in~$X$ to all the vertices in~$B_n$, as required for~$n$ by our recursion.

By construction, every two vertices $a_i,a_j$ in $A' := \{a_0, a_1,\dots\}$ have uncountably many common neighbours in~$B$: those in $B_j$ if $i<j$. By Lemma~\xxxAAfat\ applied with $B(a_i,a_j):= B_j$ for $i<j$, we deduce that $A'$~is the set of branch vertices of a fat~\TKA\ in $X$ whose subdivided edges $a_i\dots a_j$ with $i<j$ have the form $a_i b a_j$ with $b\in B_j$.%
   \COMMENT{}
   Replacing each of these paths $a_i b a_j$ with the concatenation of paths $s_i\dots v_b^i\sub S(a_i)$ and $v_b^i\dots v_b^j\sub H[V_b]$ and $v_b^j\dots s_j\sub S(a_j)$, we obtain a fat \TKA\ in~$H$ with  $s_0,s_1,\dots$ as branch vertices. (It is important here that $b$ is not just any common neighbour of $a_i$ and~$a_j$ but one in~$B_j$: only then do we know that $S(a_i)$ and $S(a_j)$ both have a leaf in~$V_b$.)
 \endproof

Let us now turn to the case of $(*)$ where $X$ is an AT-graph. As before, we shall first prove that $X$ itself contains a fat \TKA, and later refine this to a fat \TKA\ in any~$\H X$. In this second step we shall be referring to the details of the proof of the lemma below, not just to the lemma itself.

\proclaim Lemma \xxxATfat. Every AT-graph contains a fat~\TKA.

\proof
Let $X$ be an AT-graph, with Aronszajn tree $T$, say. Let us pick the branch vertices $a_0, a_1,\dots $ of our desired \TKA\ inductively, as follows.

Let $t_0$ be the root of $T_0:= T$, and $X_0 := X$. Since $X_0$ is connected,%
   \COMMENT{}
   it has a vertex $a_0$ of uncountable degree. Uncountably many of its neighbours lie above it in~$T_0$, because its down-closure is a chain and hence countable, and all its neighbours are comparable with it (by definition of an AT-graph). As levels in $T_0$ are countable, $a_0$~has a successor~$t_1$ in~$T_0$ such that uncountably many $X_0$-neighbours of $a_0$ lie above~$t_1$; let $B_0$ be some uncountable set of neighbours of~$a_0$ in~$\ucl(t_1)_{T_0}$. (We shall specify $B_0$ more precisely later.)

Let $T_1$ be the down-closure of $B_0$ in~$\ucl(t_1)_{T_0}$. Since $T_1$ is an uncountable subposet of $T_0$ with least element~$t_1$, it is again an Aronszaijn tree, and the subgraph $X_1$ it induces in $X_0$ is an AT-graph with respect to~$T_1$.

Starting with $t_0$, $T_0$ and $X_0$ as above, we may in this way select for $n=0,1,\dots$ an infinite sequence $T_0\supe T_1\supe\dots$ of Aronszajn subtrees of $T$ with roots $t_0 < t_1 <\dots$ satisfying the following:\medskip
  \item{$\bullet$} $X_n:= X[T_n]$ is an AT-graph \wrt~$T_n$;
  \item{$\bullet$} the predecessor $a_n$ of $t_{n+1}$ in~$T_n$ has an uncountable set $B_n$ of $X_n$-neighbours above $t_{n+1}$ in~$T_n$;
  \item{$\bullet$} $T_{n+1} = \ucl(t_{n+1})_{T_n} \cap \dcl(B_n)_{T_n}$.

\smallskip\noindent
By the last item above, there exists for every $b\in T_{n+1}$ a vertex $b'\in B_n\cap \ucl(b)$ (possibly $b'=b$). Applied to vertices $b$ in $B_{n+1}\sub T_{n+1}$ this means that, inductively,

\textno Whenever $i<j$, every vertex in $B_j$ has some vertex of $B_i$ in its up-closure.%
   \COMMENT{}
   &(2)

\smallskip
Let us now make $a_0, a_1,\dots$ into the branch vertices of a fat \TKA\ in~$X$. As earlier, we enumerate the desired subdivided edges as one $\omega_1$-sequence, and find independent paths $P_\alpha\sub X$ to serve as these subdivided edges recursively for all $\alpha < \omega_1$. When we come to construct the path~$P_\alpha$, beween $a_i$ and $a_j$ with $i<j$ say, we have previously constructed only the countably many paths~$P_\beta$ with $\beta<\alpha$. The down-closure $D_\alpha$ in~$T$ of all their vertices and all the $a_n$ is a countable set, since the down-closure of each vertex is a chain in~$T$ and hence countable. We can thus find a vertex $b\in B_j$ outside~$D_\alpha$, and a vertex $b'\ge b$ in~$B_i$ by~(2). The interval of $T$ between $b$ and $b'$ thus avoids~$D_\alpha$, and since it is connected in $X$ it contains the vertices of a $b'$--$b$ path $Q_\alpha$ in $X-D_\alpha$. We choose $P_\alpha := a_i b' Q_\alpha b a_j$ as the $\alpha$th subdivided edge for our fat~\TKA\ in~$X$.
 \endproof

\proclaim Lemma \xxxIATfat.
Every $\H X$, where $X$ is an AT-graph, contains a fat $\T K_{\aleph_0}$.

\proof
Let $H$ be an $\H X$ with branch sets $V_x$ for vertices $x\in X$, where $X$ is an AT-graph \wrt\ an Aronszajn tree $T$. Rather than applying Lemma~\xxxATfat\ to $X$ formally, let us re-do its proof for~$X$. We shall choose the sets~$B_n$ more carefully this time, so that we can turn the fat \TKA\ found in $X$ into one in~$H$.

Given~$n$, the set $B_n$ chosen in the proof of Lemma~\xxxATfat\ was an arbitrary uncountable set of upper neighbours of $a_n$ in~$T_n$ above some fixed successor $t_n$ of~$a_n$. We shall replace $B_n$ with a subset of itself, chosen as follows. For every $b\in B_n$, pick a vertex $v_b^n\in V_b$ that sends an edge of $H$ to a vertex~$u^n_b\in V_{a_n}$.%
   \COMMENT{}
   As in the proof of Lemma \xxxIAAfat, there is a subdivided uncountable star $S_n$ in $H$ whose leaves are among these~$v_b^n$ and all whose non-leaves, including its centre~$s_n$, lie in~$V_{a_n}$. Let us replace $B_n$ with its (uncountable) subset consisting of only those~$b$ whose $v_b^n$ is a leaf of~$S_n$.

Let $K\sub X$ be the fat \TKA\ found by the proof of Lemma~\xxxATfat\ for these revised sets~$B_n$. In order to turn $K$ into the desired \TKA\ in~$H$, we replace its branch vertices $a_n$ by the centres $s_n$ of the stars~$S_n$, and its subdivided edges $P_\alpha = a_i b' Q_\alpha b a_j$ between branch vertices $a_i, a_j$ by the concatenation of paths $s_i\dots v^i_{b'}\sub S_i$ and $Q'_\alpha = v_{b'}^i\dots v_{b}^j$ and $v^j_{b}\dots s_j\sub S_j$, where $Q'_\alpha$ is a path in $H$ expanded from~$Q_\alpha$, i.e.\ whose vertices lie in the branch sets of the vertices of~$Q_\alpha$. These paths $P'_\alpha$ are internally disjoint for distinct~$\alpha$, because the $P_\alpha$ were internally disjoint.
   \endproof

\looseproof{Proof of Theorem~\xxxMain}
Let $G$ be a connected graph without a \NST; we show that $G$ contains a fat~\TKA. By Theorem~\xxxCharacterizations, $G$~has an $X$-minor such that $X$ is either an $(\aleph_0,\aleph_1)$-graph or an Aronszajn-tree graph. Equivalently, $G$~has a subgraph $H$ that is an~$\H X$, with $X$ as above. By Lemmas~\xxxIAAfat\ and~\xxxIATfat, this subgraph~$H$, and hence~$G$, contains a fat~\TKA.
   \endproof

\beginsection References

\ref\refMetric
   R.\th Diestel, End spaces and spanning trees, \JCTB96 (2006), 846--854.

\ref\refTopSurvey
   R.\th Diestel, Locally finite graphs with ends: a topological approach. ArXiv:0912.4213 (2009).

\ref\refBook
   R.\th Diestel, {\it Graph Theory\/} (4th edition), Springer-Verlag Heidelberg, 2010.

\ref\refBGC
   R.\th Diestel \& I.B.\th Leader, A proof of the bounded graph conjecture, \Inv108 (1992), 131--162.

\ref\refNST
   R.\th Diestel \& I.B.\th Leader, Normal spanning trees, Aronszajn trees and excluded minors, \JLMS63 (2001), 16--32.

\ref\refHalin
   R.\th Halin, Simplicial decompositions of infinite graphs, in (B.\th Bollob\'as, ed): {\it Advances in Graph Theory\/}, {\sl Annals of Discrete Mathematics~\bf3} (1978).

\ref\refJung
   H.A.\th Jung, Wurzelb\"aume und unendliche {W}ege in {G}raphen, \MN41 (1969), 1--22.

\vfill
\ninepoint\noindent Version 8.4.2016
\eject\end